\documentclass{article}

\usepackage{amssymb,latexsym,amsmath}

\usepackage{graphicx}

\usepackage{graphicx,amssymb,lineno}

\begin{document}

\newcommand{\bfi}{\bfseries\itshape}

\makeatletter

\@addtoreset{figure}{section}

\def\thefigure{\thesection.\@arabic\c@figure}

\def\fps@figure{h, t}

\@addtoreset{table}{bsection}

\def\thetable{\thesection.\@arabic\c@table}

\def\fps@table{h, t}

\@addtoreset{equation}{section}

\def\theequation{\thesubsection.\arabic{equation}}

\makeatother

\newtheorem{thm}{Theorem}[section]

\newtheorem{prop}[thm]{Proposition}

\newtheorem{lema}[thm]{Lemma}

\newtheorem{cor}[thm]{Corollary}

\newtheorem{defi}[thm]{Definition}

\newtheorem{rk}[thm]{Remark}

\newtheorem{exempl}{Example}[section]

\newenvironment{exemplu}{\begin{exempl}  \em}{\hfill $\surd$

\end{exempl}}

\newcommand{\comment}[1]{\par\noindent{\raggedright\texttt{#1}

\par\marginpar{\textsc{Comment}}}}

\newcommand{\todo}[1]{\vspace{5 mm}\par \noindent \marginpar{\textsc{ToDo}}\framebox{\begin{minipage}[c]{0.95 \textwidth}

\tt #1 \end{minipage}}\vspace{5 mm}\par}

\newcommand{\ea}{\mbox{{\bf a}}}

\newcommand{\eu}{\mbox{{\bf u}}}

\newcommand{\ueu}{\underline{\eu}}

\newcommand{\ueo}{\overline{u}}

\newcommand{\oeu}{\overline{\eu}}

\newcommand{\ew}{\mbox{{\bf w}}}

\newcommand{\ef}{\mbox{{\bf f}}}

\newcommand{\eF}{\mbox{{\bf F}}}

\newcommand{\eC}{\mbox{{\bf C}}}

\newcommand{\en}{\mbox{{\bf n}}}

\newcommand{\eT}{\mbox{{\bf T}}}

\newcommand{\eL}{\mbox{{\bf L}}}

\newcommand{\eR}{\mbox{{\bf R}}}

\newcommand{\eV}{\mbox{{\bf V}}}

\newcommand{\eU}{\mbox{{\bf U}}}

\newcommand{\ev}{\mbox{{\bf v}}}

\newcommand{\eve}{\mbox{{\bf e}}}

\newcommand{\uev}{\underline{\ev}}

\newcommand{\eY}{\mbox{{\bf Y}}}

\newcommand{\eK}{\mbox{{\bf K}}}

\newcommand{\eP}{\mbox{{\bf P}}}

\newcommand{\eS}{\mbox{{\bf S}}}

\newcommand{\eJ}{\mbox{{\bf J}}}

\newcommand{\eB}{\mbox{{\bf B}}}

\newcommand{\eH}{\mbox{{\bf H}}}

\newcommand{\leb}{\mathcal{ L}^{n}}

\newcommand{\eI}{\mathcal{ I}}

\newcommand{\eE}{\mathcal{ E}}

\newcommand{\hen}{\mathcal{H}^{n-1}}

\newcommand{\eBV}{\mbox{{\bf BV}}}

\newcommand{\eA}{\mbox{{\bf A}}}

\newcommand{\eSBV}{\mbox{{\bf SBV}}}

\newcommand{\eBD}{\mbox{{\bf BD}}}

\newcommand{\eSBD}{\mbox{{\bf SBD}}}

\newcommand{\ecs}{\mbox{{\bf X}}}

\newcommand{\eg}{\mbox{{\bf g}}}

\newcommand{\paromega}{\partial \Omega}

\newcommand{\gau}{\Gamma_{u}}

\newcommand{\gaf}{\Gamma_{f}}

\newcommand{\sig}{{\bf \sigma}}

\newcommand{\gac}{\Gamma_{\mbox{{\bf c}}}}

\newcommand{\deu}{\dot{\eu}}

\newcommand{\dueu}{\underline{\deu}}

\newcommand{\dev}{\dot{\ev}}

\newcommand{\duev}{\underline{\dev}}

\newcommand{\weak}{\stackrel{w}{\approx}}

\newcommand{\mild}{\stackrel{m}{\approx}}

\newcommand{\strong}{\stackrel{s}{\approx}}

\newcommand{\weakdown}{\rightharpoondown}

\newcommand{\opg}{\stackrel{\mathfrak{g}}{\cdot}}

\newcommand{\opunu}{\stackrel{1}{\cdot}}
\newcommand{\opdoi}{\stackrel{2}{\cdot}}

\newcommand{\opn}{\stackrel{\mathfrak{n}}{\cdot}}

\newcommand{\tr}{\ \mbox{tr}}

\newcommand{\Ad}{\ \mbox{Ad}}

\newcommand{\ad}{\ \mbox{ad}}

\renewcommand{\contentsname}{ }

\title{Blurred constitutive laws  and bipotential convex covers}


\author{
 G\'ery de Saxc\'e\footnote{Laboratoire de 
  M\'ecanique de Lille, UMR CNRS 8107, Universit\'e des Sciences et 
  Technologies de Lille,
 B\^atiment Boussinesq, Cit\'e Scientifique, 59655 Villeneuve d'Ascq cedex, 
 France, e-mail: gery.desaxce@univ-lille1.fr},
 Marius Buliga\footnote{"Simion Stoilow" Institute of Mathematics of the Romanian Academy,
 PO BOX 1-764,014700 Bucharest, Romania, e-mail: Marius.Buliga@imar.ro },
Claude  Vall\'ee\footnote{Laboratoire de 
M\'ecanique des Solides, UMR 6610, UFR SFA-SP2MI, bd M. et P. Curie, 
t\'el\'eport 2, BP 30179, 86962 Futuroscope-Chasseneuil cedex, 
France, e-mail: vallee@lms.univ-poitiers.fr}
}

\date{This version: April 22, 2009}


\maketitle

{\bf MSC-class:} 49J53; 49J52; 26B25
\begin{abstract}

In many practical situations, incertitudes affect the mechanical behaviour 
that is given by a family of graphs instead of a single one. 
In this paper, we show how the bipotential method is able to capture such blurred 
constitutive laws, using  bipotential convex covers. 

\end{abstract}

\section{Introduction}

The constitutive laws of the materials can be represented, as in Elasticity, 
by a univalued mapping $T:X\rightarrow Y$ or, as in in Plasticity, can be 
put in the form of a multivalued operator $T:X\rightarrow 2^{Y}$.  
Equivalently, a constitutive law can be seen as the graph 
$M \subset X \times Y$ of the operator $T$, that is  
$\displaystyle M = \, Graph(T) = \left\{ (x,y) \mbox{ : } y \in T(x) \right\}$. 
$X$ and $Y$ are spaces of dual variables, for example $X$ may be a space of 
stresses and $Y$ may be a space of deformation rates. The duality between these 
spaces is a function $\langle \cdot, \rangle : X\times Y \rightarrow
\mathbb{R}$.

If the graph $M$ is maximal cyclically monotone, 
then there is a convex and lower semi-continuous (l.s.c.) 
function $\phi:X\rightarrow \bar{\mathbb{R}}$, called a 
\textbf{superpotential} (or pseudo-potential), such that  $\displaystyle 
M = \, Graph(\partial \phi)$, where  $\partial \phi$ is 
the  subdifferential of $\phi$ and the function $\phi$ is determined by the
graph $M$, up to an additive constant.  

The constitutive laws of dissipative materials admitting a superpotential can be
put into the form $y \in \partial \phi(x)$. Such constitutive laws 
are often qualified as standard  \cite{Halp JM 75} and the law is said to 
be a normality law, a subnormality law or an associated law. 
However, many experimental laws proposed these last decades, 
particularly in Plasticity, are non associated. 
For such laws, we proposed in \cite{saxfeng} a suitable modelization 
thanks to a function called bipotential. 

The laws admitting a bipotential are called laws of implicit standard 
materials because they have the form $y \in \partial b(\cdot , y)(x)$ , which 
is a subnormality law but the relation between $x$ and $y$ is implicit.

The bipotential theory allows, 
in connection with the calculus of variation, 
to model a wide spectrum of non associated constitutive laws.  
 Examples of such  non associated constitutive laws are:
non-associated Dr\"ucker-Prager \cite{sax boush KIELCE 93}  and 
Cam-Clay models \cite{sax BOSTON 95} in soil mechanics, 
cyclic Plasticity (\cite{sax CRAS 92},\cite{bodo sax EJM 01}) 
and Viscoplasticity \cite{hjiaj bodo CRAS 00} of metals with non linear 
kinematical hardening rule, Lemaitre's damage law \cite{bodo}, the coaxial 
laws (\cite{sax boussh 2},\cite{vall leri CONST 05}), the Coulomb's friction law \cite{saxfeng}, 
\cite{sax CRAS 92}, 
\cite{boush chaa IJMS 02}, \cite{feng hjiaj CM 06}, \cite{fort hjiaj CG 02}, 
\cite{hjiaj feng IJNME 04}, \cite{sax boush KIELCE 93}, \cite{sax feng IJMCM 98}, \cite{laborde}.
A complete survey can be found in \cite{sax boussh 2}.  In the previous works, robust numerical algorithms were proposed to solve 
structural mechanics problems.

The cornerstone inequality in the definition of bipotentials (definition
\ref{def2} (b)) extends the  Fenchel's inequality. In particular, to any 
superpotential $\phi$ is associated  the 
\textbf{separable bipotential}:
\begin{equation}
   b(x,y) = \phi(x) + \phi^{*}(y)
\label{separable bipotential} 
\end{equation} 
where $\phi^{*}$ is the Fenchel conjugate of $\phi$ (with respect to the duality
between the spaces $X$ and $Y$). 
The implicit 
subnormality law $y \in \partial b(\cdot , y)(x)$ becomes the associated 
law $y \in \partial
\phi(x)$. However, there are many bipotentials which can not be expressed in the
form (\ref{separable bipotential}).

For all the particular constitutive laws previously mentioned, 
the bipotentials were heuristically constructed, without knowing beforehand 
the conditions under which ones the law admits a bipotential, nor a systematic 
algorithm to construct this bipotential. 
 
In  \cite{bipo1} we solved two key problems: (a) when the graph of a given multivalued operator can be expressed as the set of critical points of a bipotentials, and (b) a method of construction of a bipotential associated 
(in the sense of point (a)) to a multivalued, typically non monotone, operator. The main tool was the notion of {\bf convex lagrangian cover} of the graph of the multivalued operator, and a related notion of implicit convexity of this cover. 
The results of \cite{bipo1}  apply  only to bi-convex, bi-closed graphs 
(for short BB-graphs) admitting  at least one convex lagrangian cover by 
{\bf maximal cyclically monotone graphs}. 
This is a rather large class of graph of multivalued operators but important 
applications to the mechanics, such as the bipotential associated to contact 
with friction \cite{saxfeng}, are not in this class. 

In more recent papers \cite{bipo2}, \cite{bipo3}, we  proposed an extension of  the method 
presented in \cite{bipo1} to a more general class of BB-graphs. This is done 
in two steps. In the first step we proved that the intersection of two maximal cyclically
monotone graphs  is the critical set of a bipotential if and only if a condition formulated in terms of 
the inf convolution of a family of convex lsc functions is true \cite{bipo3} . 
In the second step we extended  the main result of \cite{bipo1} by replacing 
the notion of convex lagrangian cover with the one of 
{\bf bipotential convex cover} (definition \ref{defcover}). In this way we were able 
to apply our results to the bipotential for the Coulomb's friction law.

The purpose of this paper is to describe a new application of bipotentials. 
In many practical situations, incertitudes affect the mechanical behaviour. 
In other words, we tolerate indeterminacy of the constitutive law which is 
represented by a family of graphs instead of a single one. 
In particular, when no solution can be found for ill-posed problems, 
relaxation of stronger conditions on the material behavior would allow to 
provide at least an approximate solution. 
Our aim now is to show how the bipotential is able to capture such blurred 
constitutive law, by using  bipotential convex covers.

The main results of the paper are Propositions \ref{pelas}, \ref{pplas} and 
\ref{pfric}. In Proposition \ref{pelas} we find a bipotential for a blurred
Elasticity law, or equivalently, we show that such law can be expressed as a 
implicit subnormality law. We pass then to a more difficult blurred 
Plasticity law, with a variable yielding threshold $\eta$ taking arbitrary 
values in  $[\lambda_{-}, \lambda_{+}]$. For this law  we achieve the same as previously in Proposition \ref{pplas}. 

The third result concerns a blurred Coulomb friction law. The law of 
unilateral contact with Coulomb's dry friction  is a typical example of  a 
non associated constitutive law in mechanics which admits a bipotential 
\cite{saxfeng}. It is to be remarked that the Coulomb friction law does not 
have a separated bipotential, because the graph of the law is not even monotone.
In Proposition \ref{pfric} we are able to provide a bipotential formulation 
for  a blurred Coulomb friction
law with  arbitrary values for the friction coefficient $\mu$ 
 in a range $\left[ \mu_{-}, \mu_{+} \right]$.

\paragraph{Aknowledgements.} The first two authors acknowledge the support 
from the European Associated Laboratory "Math Mode" associating the 
Laboratoire de Math\'ematiques de l'Universit\'e Paris-Sud (UMR 8628) and the 
"Simion Stoilow" Institute of Mathematics of the Romanian Academy. 
The first author thanks to M. Jean and M. Raous of the Laboratoire de 
M\'ecanique et d'Acoustique (UPR CNRS, Marseille), the first one for the idea at 
the origin of this work and the second one for the kind invitation to give a 
seminar which led to discussions where  this idea was suggested.

\section{Bipotentials}

$X$ and $Y$ are topological, locally convex, real vector spaces of dual 
variables $x \in X$ and $y \in Y$, with the duality product 
$\langle \cdot , \cdot \rangle : X \times Y \rightarrow \mathbb{R}$. 
We shall suppose that $X, Y$ have topologies compatible with the duality 
product, that is: any  continuous linear functional on $X$ (resp. $Y$) 
has the form $x \mapsto \langle x,y\rangle$, for some $y \in Y$ (resp. 
$y \mapsto \langle x,y\rangle$, for some  $x \in X$). We use the notation: $\displaystyle \bar{\mathbb{R}} = \mathbb{R}\cup \left\{ +\infty \right\}$. For any convex and closed set $A \subset X$, its  indicator function,  
$\displaystyle \chi_{A}$, is defined by 
$$\chi_{A} (x) = \left\{ \begin{array}{ll}
0 & \mbox{ if } x \in A \\ 
+\infty & \mbox{ otherwise } 
\end{array} \right. $$
The subgradient of a function $\displaystyle \phi: X \rightarrow
\bar{\mathbb{R}}$ at a point $x \in X$ is the (possibly empty) set: 
$$\partial \phi(x) = \left\{ u \in Y \mid \forall z \in X  \  \langle z-x, u \rangle \leq \phi(z) - \phi(x) \right\} \  .$$

\begin{defi} A {\bf bipotential} is a function $b: X \times Y \rightarrow
 \bar{\mathbb{R}}$, with the properties: 
\begin{enumerate}
\item[(a)] $b$ is convex and lower semicontinuous in each argument; 
\item[(b)] for any $x \in X , y\in Y$ we have $\displaystyle b(x,y) \geq \langle x, y \rangle$; 
\item[(c)]  for any $(x,y) \in X \times Y$ we have the equivalences: 
\begin{equation}
y \in \partial b(\cdot , y)(x) \ \Longleftrightarrow \ x \in \partial b(x, \cdot)(y)  \ \Longleftrightarrow \ b(x,y) = 
\langle x , y \rangle \ .
\label{equiva}
\end{equation}
\end{enumerate}
The {\bf graph} of $b$ is 
\begin{equation}
M(b) \ = \ \left\{ (x,y) \in X \times Y \ \mid \ b(x,y) = \langle x, y \rangle \right\} \  .
\label{mb}
\end{equation}
\label{def2}
\end{defi}

If the graph $M$ of a law is the graph of a bipotential $b$, we say that the 
law (the graph) admits a bipotential.

For any graph $M \subset X \times Y$, we can introduce the sections 
$\displaystyle M (x) \ = \ \left\{ y \in Y \mid (x,y) \in M \right\}$ and 
$ M^{*}(y) \ = \ \left\{ x \in X \mid (x,y) \in M \right\}$. Hence the operator  
$T$ assigns to each $x \in X$ the section $M (x)$ and the inverse law assigns 
to each $y \in Y$ the section $M^{*} (y)$. 

Let a constitutive law be given by a graph $M$. Does it admit a bipotential? 
The existence problem is easily settled by the following result. 

\begin{thm}
 Given  a non empty set $M \subset X \times Y$, there is a bipotential $b$ 
such that $M=M(b)$ if and only if for any $x \in X$ and $y \in Y$ the sections 
$M(x)$ and $M^{*}(y)$ are convex and closed. 
 \label{thm1}
 \end{thm}
 
 The proof can be found in \cite{bipo1}. Then we say that $M$ is bi-convex and 
 bi-closed, or in short that $M$ is a {\bf BB-graph}. 
 
 If the law is represented by a BB-graph, then a closely related topic  
 is to know whether the bipotential is unique. The answer is no. 
 The proof of the previous result is based on the introduction of the 
 bipotential
$$ b_{\infty} (x,y) = \left\langle x,y \right\rangle  + \chi_M (x,y) .$$
If the graph  $M = Graph(\partial \phi)$ is cyclically monotone maximal, then 
it admits at least two distinct bipotentials: the separable bipotential 
defined by \ref{separable bipotential} and $b_{\infty} $. 
Therefore the graph of the law alone is not  sufficient to uniquely define the 
bipotential.  

\section{Bipotential convex covers}

 Theorem \ref{thm1} does not give a 
satisfying bipotential for a given multivalued constitutive law, 
because the bipotential $b_{\infty}$ is somehow degenerate. 
We would like to find a method of construction of  bipotentials  
which for a given BB-graph $M$ will return a bipotential 
 $b$ which is not everywhere infinite outside the graph $M$, and such that if 
 $M$ is maximal cyclically monotone then the method will give us a separable
 bipotential. 
 
 We saw that the graph alone is not sufficient to construct interesting 
 bipotentials. We need more information to start from. This is provided by the 
 notion of bipotential convex cover. 

Let $Bp(X,Y)$ be the set of all bipotentials $b: X \times Y \rightarrow
\bar{\mathbb{R}}$. We shall need the following definitions (4.1 and 4.1
\cite{bipo3}).

\begin{defi}
Let $\Lambda$ be an arbitrary non empty set and $V$ a real vector space. The 
function $f:\Lambda\times V \rightarrow \bar{\mathbb{R}}$ is 
{\bf implicitly  convex} if for any two elements 
$\displaystyle (\lambda_{1}, z_{1}) , 
(\lambda_{2},  z_{2}) \in \Lambda \times V$ and for any two numbers 
$\alpha, \beta \in [0,1]$ with $\alpha + \beta = 1$ there exists 
$\lambda  \in \Lambda$ such that 
\begin{equation}
   f(\lambda, \alpha z_{1} + \beta z_{2}) \ \leq \ \alpha 
   f(\lambda_{1}, z_{1}) + \beta f(\lambda_{2}, z_{2}) \quad .
\label{ineqimpl} 
\end{equation} 
\label{defimpl}
\end{defi}

\begin{defi}   A {\bf bipotential
convex cover} of the non empty set $M$ is a function   
$\displaystyle \lambda \in \Lambda \mapsto b_{\lambda}$ from  $\Lambda$ with 
values in the set  $Bp(X,Y)$, with the 
properties:
\begin{enumerate}
\item[(a)] The set $\Lambda$ is a non empty compact topological space, 
\item[(b)] Let $f: \Lambda \times X \times Y \rightarrow \mathbb{R} \cup
\left\{ + \infty \right\}$ be the function defined by 

$$f(\lambda, x, y) \ = \ b_{\lambda}(x,y) .$$

Then for any $x \in X$ and for any $y \in Y$ the functions 
$f(\cdot, x, \cdot): \Lambda \times Y \rightarrow \bar{\mathbb{R}}$ and 
$f(\cdot, \cdot , y): \Lambda \times X \rightarrow \bar{\mathbb{R}}$ are  lower 
semi continuous  on the product spaces   $\Lambda \times Y$ and respectively 
$\Lambda \times X$ endowed with the standard topology, 
\item[(c)] We have $\displaystyle M  \ = \  \bigcup_{\lambda \in \Lambda} 
M(b_{\lambda})$.
\item[(d)] with the notations from point (b), the functions $f(\cdot, x, \cdot)$ 
and $f(\cdot, \cdot , y)$ are implicitly convex in the sense of Definition 
\ref{defimpl}.
\end{enumerate}
\label{defcover}
\end{defi}

The next theorem,  \cite{bipo3} theorem 4.6, is the key result needed further.

\begin{thm} Let $\displaystyle \lambda \mapsto b_{\lambda}$ be a bipotential 
convex cover  of 
 the graph $M$ and $b: X \times Y \rightarrow R$ defined by
\begin{equation}
b(x,y) \ = \ \inf \left\{ b_{\lambda}(x,y) \ \mid \  \lambda \in \Lambda \right\} \ . 
\end{equation}
Then $b$ is a bipotential and $M=M(b)$. 
\label{thm2}
\end{thm}

The result is rather surprising because an inferior envelop of functions, 
even convex, is not generally a convex function. The property (d)  of the 
Definition \ref{defcover} is essential to ensure the convexity properties of 
$b$.

\section{Application to a Elasticity law with thick line}

The finite element method is a numerical method of discretizing continua, 
widely used today for solving structural mechanics problems. 
The accuracy of the approximative solution can be controlled by computing a 
posteriori estimators. 
Three different approaches were proposed. Contrary to the methods based on 
the equilibrium residuals (\cite{ER1}, \cite{ER2}, \cite{ER3}, \cite{ER4}, 
\cite{ER5}), the ones using smoothing techniques (\cite{ST1}, \cite{ST2}, 
\cite{ST3}) and the dual analysis based on upper and lower bounds for the 
energy (\cite{DA1}, \cite{DA2}), the method of the error on the constitutive 
law  (or constitutive relation error) is based on mechanical concepts and can 
be more naturally extended to non linear problems of evolution 
(\cite{LAD1}, \cite{LAD2}, \cite{LAD3}, \cite{LAD4}, \cite{LAD5}, 
\cite{LAD6}, \cite{LAD7}, \cite{LAD8}). 
We use such an idea in the sense that we admit an error on the elastic 
law instead of satisfying it exactly. We shall use bipotentials constructed 
from bipotential convex covers, in order to formulate the elasticity law with
thick line as a implicit subnormality law. 

We shall take $X=Y=\mathbb{R}^{n}$ and the duality product is the usual scalar product in $\mathbb{R}^{n}$. let us consider the elastic linear law $y = \lambda x$ with $\lambda > 0$ which is the most simple example of linear elastic law where the dual variables $x$ and $y$ are vectors. In the present application, the material parameter, the "elastic modulus" $\lambda$, has a fixed value but we allow some error on the constitutive law $a = y - \lambda x$ with a fixed margin of tolerance $\epsilon > 0$. 
In other words, we considered the blurred elastic law described by the BB-graph:
\begin{equation}
    M = \left\{ (x,y)\in X \times Y \   \mid \   \parallel y - \lambda x  \parallel  \leq \epsilon \right\}\ . 
\label{blurred elastic law} 
\end{equation} 



\begin{figure}[h!]
\centering
 \includegraphics[width=0.9\textwidth]{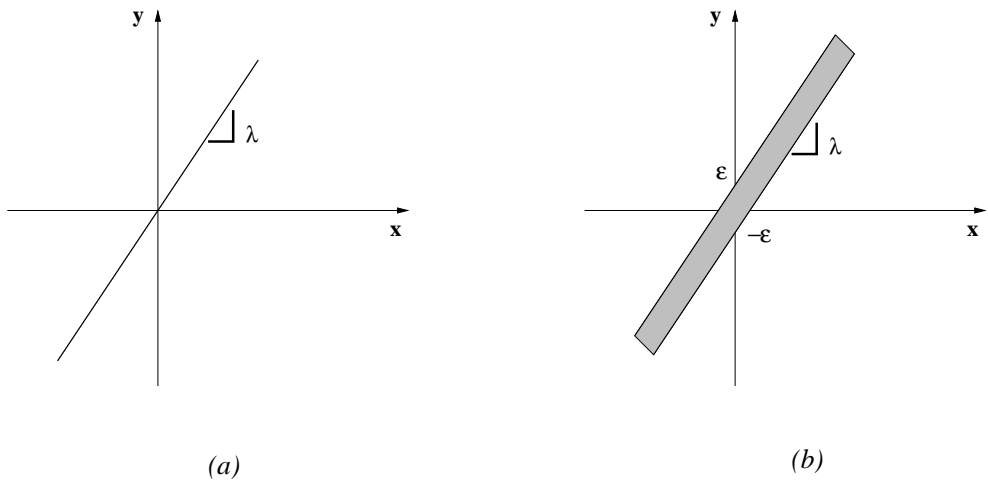}
 \caption{Elasticity law. (a): ideal law with thin line. (b): law with thick line}
\label{fig 1}
\end{figure}


In figure \ref{fig 1}, a pictural representation of the blurred law is a graph 
with "thick line" (displayed at the right) in contrast to the ideal law of the 
previous section with a "thin line" (displayed at the left). 

\begin{prop}
The blurred elastic law with the graph $M$ (\ref{blurred elastic law}) is
represented by the bipotential: 
\begin{equation}
b (x,y) = \left\langle  x,y \right\rangle  
                     + \frac{1}{2 \lambda}\ \left(\  \left( \parallel y - \lambda x \parallel - \epsilon\right)_+ \right)^2
\label{bipoelas}
\end{equation}
with the notation $\displaystyle \left(\alpha\right)_{+} = \max(\alpha, 0)$. 

Equivalently, the relation  $\displaystyle \|y - \lambda x \| \leq \varepsilon$ 
can be put in the form $y \in \partial b(\cdot,y)(x)$, with $b$ the bipotential 
defined by (\ref{bipoelas}).
\label{pelas}
\end{prop}

\paragraph{Proof.} The graph $M$ is clearly bi-convex and bi-closed. 
We construct a bipotential convex cover of $M$
 by assigning  to the parameter  $a \in B (\epsilon) $ the set
$$ M_a = \left\{ (x,y)\in X \times Y \   \mid \  y = \lambda x   \right\}\ ,  $$
which can be seen as the graph of the elastic law $y = \lambda x + a$ with an "initial stress" $a$. it is clear that $M_a$ is the subdifferential $\partial \phi_a$ of the potential:
$$ \phi_a (x) = \frac{\lambda}{2} \parallel x \parallel^2 + \left\langle x,a \right\rangle\ . $$
Its Fenchel conjugate is:
$$ \phi_{a}^{*} (y) = \frac{1}{2 \lambda} \parallel y - a \parallel^2\ . $$
 Let  $\left\{ b_{a} \mbox{ : } a \in B (\epsilon) \right\}$ be the collection of the separated bipotentials:
$$ b_{a}(x,y) \ = \phi_a (x) + \phi_{a}^{*} (y)
                        \ = \ \frac{\lambda}{2} \parallel x \parallel^2 + \left\langle x,a \right\rangle
                  +  \frac{1}{2 \lambda} \parallel y - a \parallel^2\ ,$$
$$ b_{a}(x,y) \ = \ \left\langle x,y \right\rangle 
                               + \frac{1}{2 \lambda} \parallel y - a - \lambda x \parallel^2\ $$
We want to verify that it defines a bipotential convex cover of $M$. The 
conditions (a) to (c) of Definition \ref{defcover} are fulfilled. 
We have to prove the last condition (d). 
For the implicit convexity of $f(\cdot, x, \cdot)$, the inequality 
\ref{ineqimpl} is then:  for any $x \in X$, $y_1,y_2 \in Y$, 
$a_1, a_2 \in B (\epsilon)$, $\alpha, \beta \in [0,1]$, $\alpha + \beta = 1$,
there exists $a \in B(\varepsilon)$ such that: 
$$ \alpha \parallel y_1 - a_1 - \lambda x \parallel^2 + \beta \parallel y_2 - a_2 - \lambda x \parallel^2
        \geq \parallel \alpha y_1 + \beta y_2 - a - \lambda x \parallel^2\ .$$
Because the square of the norm is convex, an obvious choice for $a$ is 
$a = \alpha a_1 + \beta a_2$.  
For the implicit convexity of $f(\cdot, \cdot, y)$, the demonstration is similar.
We apply  Theorem \ref{thm2} and we obtain a  bipotential of the graph $M$ with
the form:
$$b(x,y)  \ = \ \inf \left\{ b_{a}(x,y) \ \mbox{ : } a \in B (\epsilon) \right\}\
                  = \ \inf \left\{ \frac{\lambda}{2} \parallel x \parallel^2 + \left\langle x,a \right\rangle
                  +  \frac{1}{2 \lambda} \parallel y - a \parallel^2
 \mbox{ : } \parallel a \parallel \leq \epsilon) \right\} \ .$$
We shall prove now that $b$ has the form (\ref{bipoelas}). 
Two events have to be considered:
\begin{enumerate}
 \item[-] The infimum is realized at $a$ such that 
 $\parallel a \parallel < \epsilon $. Hence $b_{a}(x,y)$ is stationary with 
 respect to $a$, that is: $ a = y - \lambda x$. 
Eliminating $a$ by this relation, a straightforward calculation shows that the 
infimum is:
$$ b(x,y) = \left\langle   x,y \right\rangle\ .$$
 \item[-] Otherwise, introducing a Lagrange multiplier $\eta$, we have:
$$ b (x,y) \ =  \  \sup \left\{ \inf 
                                                 \left\{ \ \ b_{a}(x,y) + \eta\ \left( \parallel a \parallel^2 - \epsilon^2 \right) 
                                                           \ \mbox{ : } a \in Y
                                            \right\}
                                           \ \mbox{ : } \eta \geq 0 
                                  \right\}  $$
\end{enumerate}
The stationarity condition with respect to $a$ gives:
\begin{equation}
   a = \frac{y - \lambda x}{1 + 2 \eta}\ .
\label{a} 
\end{equation}   
Hence the constraint $ \parallel a \parallel = \epsilon$ allows to deduce the value of the Lagrange multiplier:
\begin{equation}
   \eta = \frac{1}{2} \left( \frac{1}{\epsilon} \parallel y - \lambda x \parallel - 1 \right) \ .
\label{eta} 
\end{equation}
Introducing expression \ref{eta} into \ref{a} leads to:
$\displaystyle  a = \epsilon \frac{y - \lambda x}{\parallel y - \lambda x
\parallel}$. 
Eliminating $a$ by this relation gives the value of the infimum:
$$ b (x,y) =  \left\langle  x,y \right\rangle  
                     + \frac{1}{2 \lambda}\ \left(  \parallel y - \lambda x \parallel - \epsilon \right)^2
 \ .$$
We verify immediately that $b$ has indeed the expression (\ref{bipoelas}. \quad
$\square$

The same reasoning may be performed for the more realistic elasticity law 
$y = K x$ with $K$ a elasticity tensor, but the computations are more involved.

\section{Application to a blurred Plasticity law}

We want now to extend the previous ideas to non smooth constitutive laws. 
Let us consider  the Plasticity law with a yielding threshold $\lambda$ for 
which the plastic domain is the the closed ball $B(\lambda)$ of center $0$ and 
radius $\lambda$: $\displaystyle B(\lambda) \ = \ \left\{ y \in Y \mbox{ : } 
\| y \| \leq \lambda \right\}$. 
The Plasticity law is given by the maximal cyclically monotone graph:
$$  M_{\lambda}  \ = \ \left\{ (0,y) \in X \times Y \mbox{ : } \| y \| < \lambda \right\}
                                     \ \ \cup\ \  \left\{ (x,y)\in X \times Y \  \mbox{ : }  \  
                       \parallel y \parallel = \lambda, \exists \eta \geq 0, x = \eta y  \right\}   \ . $$
The graph $M_{\lambda}$ is the subdifferential $\partial \phi_{\lambda}$ of 
the potential:
\begin{equation}
   \phi_{\lambda} (x) = \  \ \lambda \| x\| \ .
\label{power of dissipation} 
\end{equation}  
Its Fenchel conjugate is: $\displaystyle  \phi_{\lambda}^{*} (y) = \chi_{B(\lambda)}(y)$. 
We allow some error on the constitutive law:
$$ y \in \partial \phi_{\lambda} (x) + a$$ 
with a fixed margin of tolerance $\epsilon > 0$ on the norm of $a$. In other words, we consider the 
blurred plastic law described by the graph:
\begin{equation}
    M = \left\{ (x,y)\in X \times Y \  \mid \  \exists\ a \in B (\epsilon),\ y
    \in \partial \phi_{\lambda} (x) + a  \ \right\}\ . 
\label{blurred plastic law} 
\end{equation} 


\begin{figure}[h!]
 \centering
 \includegraphics[width=0.9\textwidth]{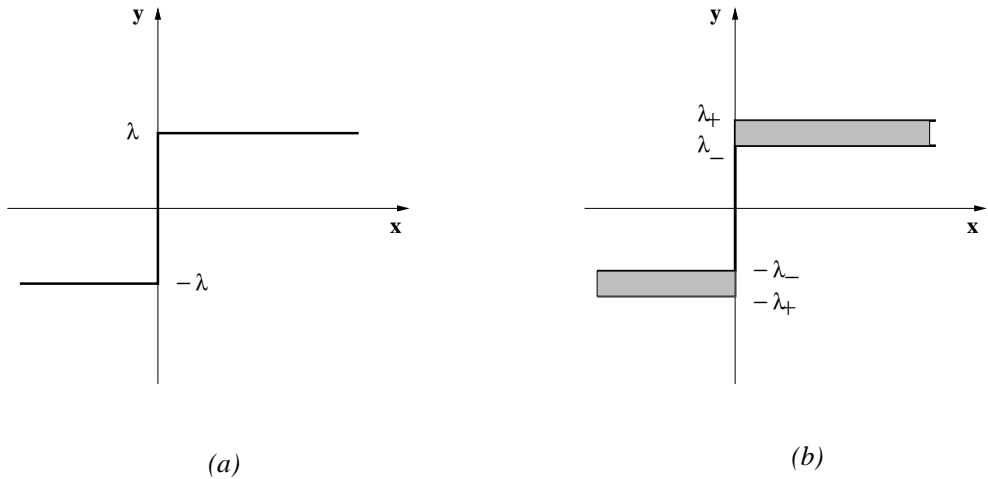}
 \caption{Plasticity law. (a): ideal law with thin line. (b): law with thick line}
\label{fig 2}
\end{figure}


In figure \ref{fig 2}, a pictural representation of the blurred law is a graph with "thick line" 
(displayed at the right) in contrast to the ideal law  with a "thin line" 
(displayed at the left).

\begin{prop}
The blurred plastic law with the graph $M$ (\ref{blurred plastic law}) is
represented by the bipotential: 
\begin{equation}
b (x,y) = sup \left( \lambda_{-} , \parallel y \parallel \right)  \parallel x \parallel + \chi_{B(\lambda_{+})}(y)
\label{bipoplas}
\end{equation}
with the notation $\lambda_{\pm} =  \lambda \pm \epsilon $. 

Equivalently, the relation  $ y \in \partial \phi (x) + a$, $\|a\| \leq
\varepsilon$,  
can be put in the form $y \in \partial b(\cdot,y)(x)$, with $b$ the bipotential 
defined by (\ref{bipoplas}).
\label{pplas}
\end{prop}

\paragraph{Proof.} 
As previously, we intend to construct a bipotential for $M$ from a bipotential
convex cover. This time  $a$ is not the good 
parameter for a bipotential convex cover because it is difficult 
to check  all the conditions of 
Theorem \ref{thm2} for this parameter. 

The bipotential convex cover will have as parameter a  
variable yielding threshold $\eta$ taking arbitrary values in  
$[\lambda_{-}, \lambda_{+}]  $ where $\lambda_{\pm} =  \lambda \pm \epsilon $. 
The graph $M$ (\ref{blurred plastic law} ) admits the description:
\begin{equation}
M = \left\{ (0,y)\in X \times Y \   \mid \    \parallel y \parallel <  \lambda_{-} \right\} \cup 
\label{plasticity law with thick line}
\end{equation} 
$$ \cup \left\{ (x,y)\in X \times Y \   \mid \   \lambda_{-} \leq \parallel y \parallel \leq \lambda_{+}\
                       and \ \exists \eta \geq 0, \ x = \eta y \right\} $$

 Let  $\left\{ b_{\eta} \mbox{ : } \eta \in [\lambda_{-}, \lambda_{+}] 
 \right\}$ be the collection of  separated bipotentials:
$$ f(\eta, x, y) = b_{\eta}(x,y) \ = \ \eta \| x \|+ \chi_{B(\eta)}(y)\ .$$
We shall verify  that this is  a bipotential convex cover of $M$. 
Indeed, the conditions (a) to (c) of Definition \ref{defcover} are obviously 
fulfilled. With $f (\eta,x,y) = b_\eta (x,y) $, we have to prove the last 
condition (d). For the implicit convexity of $f(\cdot, \cdot , y)$ we have to
prove that for any $x_1,x_2 \in X$, $y \in Y$, $\eta_{1}, \eta_{2} \in [\lambda_{-}, \lambda_{+}]$, 
$\alpha, \beta \in [0,1]$, $\alpha + \beta = 1$ there exists $\eta \in
[\lambda_{-}, \lambda_{+}]$ such that: 
 $$\alpha \eta_1 \parallel x_1 \parallel + \beta \eta_2 \parallel x_2 \parallel 
   + \chi_{B(\eta_1)}(y)+ \chi_{B(\eta_2)}(y) \geq
   \eta \parallel \alpha x_1 + \beta x_2 \parallel + \chi_{B(\eta)}(y)\ .$$
We choose $ \eta = \min ( \eta_1, \eta_2)$. The previous inequality 
becomes: for any  $\parallel y \parallel  \leq \eta $   
$$ \alpha \eta_1 \parallel x_1 \parallel + \beta \eta_2 \parallel x_2 \parallel  \geq
   \eta \parallel \alpha x_1 + \beta x_2 \parallel\ ,$$
which is true because of  the convexity of the norm. 
 
For the implicit convexity of $f(\cdot, x, \cdot)$ we have to prove that:  
for any $x \in X$, $y_1,y_2 \in Y$,  $\eta_{1}, \eta_{2} \in 
[\lambda_{-}, \lambda_{+}]$, $\alpha, \beta \in [0,1]$, $\alpha + \beta = 1$there exists $\eta \in
[\lambda_{-}, \lambda_{+}]$ such that:
$$ (\alpha \eta_1 + \beta \eta_2)\  \parallel x \parallel 
      + \chi_{B(\eta_1)}(y_1)+ \chi_{B(\eta_2)}(y_2) \geq
      \eta \parallel x \parallel + \chi_{B(\eta)}(\alpha y_1 + \beta y_2)$$
This time we choose $\eta = \alpha \eta_1 + \beta \eta_2$. 
Indeed,  for any  $\parallel y_1 \parallel  \leq \lambda_1 $ and 
$\parallel y_2 \parallel  \leq \lambda_2 $, we have:
$$\parallel \alpha y_1 + \beta y_2 \parallel
    \  \leq\  \alpha \parallel y_1 \parallel +  \beta \parallel y_2 \parallel 
     \ \leq\  \alpha \eta_1 + \beta \eta_2 = \eta$$
therefore the inequality we want to prove becomes trivial: 
$$ (\alpha \eta_1 + \beta \eta_2)\  \parallel x \parallel  \geq
   \eta \parallel x \parallel \ .$$

Hence we are in the conditions of applying Theorem \ref{thm2}. 
The graph $M$ admits the bipotential 
$$b(x,y)  \ = \ \inf \left\{ \lambda \|x\| + \chi_{B(\lambda)}(y) 
\mbox{ : } \lambda \in [\lambda_{-}, \lambda_{+}] \right\} \ .$$
In order to compute $b$, three events have to be considered:
\begin{enumerate}
 \item[-] for $\parallel y \parallel >  \lambda_{+}$, $ b (x,y) \ = +\infty\ $ , 
 \item[-] for $\lambda_{-} \leq \parallel y \parallel \leq  \lambda_{+}$, $b(x,y)  \ = \ \inf \left\{ \lambda \|x\|  \mbox{ : }  \parallel y \parallel < \lambda \leq  \lambda_{+} \right\}  \ = \ \|y\| \|x\| \ $, 
 \item[-] for $\parallel y \parallel \leq  \lambda_{-}$, $b(x,y)  \ = \ \inf \left\{ \lambda \|x\|  \mbox{ : }    \lambda_{-} \leq \lambda \leq  \lambda_{+} \right\}  \ = \  \lambda_{-} \|x\| \ $.
\end{enumerate}
The proof is done. \quad $\square$

\section{Application to a blurred Coulomb's friction law}

The law of unilateral contact with Coulomb's dry friction  is a typical 
example of what is called a non associated constitutive law in mechanics. 
Despite of its rather complex structure, it is worthwhile to have interest 
in it because of its importance in many practical problems. 

We shall not discuss here the phenomenal and experimental aspects but only 
the mathematical modeling with respect to the bipotential theory. 
To be short, the space $X = \mathbb{R}^{3}$ is the one of relative velocities 
between points of two bodies, and the space $Y$, identified also to 
$\mathbb{R}^{3}$, is the one of the contact reaction stresses. The duality 
product is the usual scalar  product. We put
$$(x_{n},x_{t})\in X = \mathbb{R} \times \mathbb{R}^{2}, \quad 
    (y_{n},y_{t})\in Y = \mathbb{R} \times \mathbb{R}^{2}\ ,$$
where $x_{n}$ is the gap velocity, $x_{t}$ is the sliding velocity, 
$y_{n}$ is the contact pressure and $y_{t}$ is minus the friction stress. 
The friction coefficient is $\mu > 0$. The graph of the law of unilateral 
contact with Coulomb's dry friction is defined as the union of three sets, 
respectively corresponding to the 'body separation', the 'sticking' and the 
'sliding'. 
\begin{equation}
M_{\mu} = \left\{ (x,0)\in X \times Y \   \mid \ x_{n} < 0 \right\} 
     \cup \left\{ (0,y)\in X \times Y \   \mid \  
                       \parallel y_{t} \parallel \leq \mu y_{n}  \right\} 
		       \cup
\label{Coulomb friction contact law}
\end{equation} 
$$\cup \left\{ (x,y) \in X \times Y \   \mid \  x_{n} = 0, \ x_{t} \neq 0, \ 
                       y_{t} = \mu y_{n} \dfrac{x_{t}}{\parallel x_{t}
		       \parallel}   \right\} $$
It is well known that this graph is not monotone, then not cyclically 
monotone. As usual, we introduce Coulomb's cone
$$ K_{\mu} =  \left\{ (y_{n},y_{t})\in Y \   \mid \  
                      \parallel y_{t} \parallel \leq \mu y_{n}  \right\}  ,$$
and its conjugate cone
$$ K_{\mu}^* =  \left\{ (x_{n},x_{t})\in X \   \mid \  
                      \mu \parallel x_{t} \parallel  + x_{n} \leq 0  \right\}  .$$
In particular, we have
$$ K_{0} =  \left\{ (y_{n},0)\in Y \   \mid \  y_{n} \geq 0  \right\}  , 
\quad 
     K_{0}^* =  \left\{ (x_{n},x_{t})\in X \   \mid \   x_{n} \leq 0  \right\}  .$$
In \cite{bipo2}, \cite{bipo3}, we obtained  as application of Theorem 
\ref{thm2} the following expression: 
\begin{equation}
   b_{\mu} (x,y) \ =  \mu y_{n} \parallel x_{t} \parallel + \chi_{K_{\mu}} (y) +  \chi_{ K^{*}_{0} } (x) 
\label{Coulomb bipotential} 
\end{equation} 
recovering the bipotential firstly obtained in heuristic way in \cite{saxfeng}.

 We shall  modify the unilateral contact law with Coulomb's dry 
 friction by allowing arbitrary values for the friction coefficient $\mu$ 
 in a range $\left[ \mu_{-}, \mu_{+} \right] $, that leads to the blurred 
 friction law represented by the graph:
\begin{equation}
M = \left\{ (x,0)\in X \times Y \   \mid \ x_{n} < 0 \right\} \cup
\label{blurred Coulomb friction contact law}
\end{equation} 
$$  \cup \left\{ (0,y)\in X \times Y \   \mid \  
                      \mu_{-} y_{n} \leq   \parallel y_{t} \parallel \leq \mu_{+} y_{n},
                       \ \  \exists \eta \geq 0,\ \ x_t = \eta y_t \right\} $$

\begin{prop}
The blurred plastic law with the graph $M$ (\ref{blurred Coulomb friction contact law}) is
represented by the bipotential: 
\begin{equation}
b (x,y) = \sup \left( \mu_{-} y_n , \parallel y_t \parallel \right)  \parallel x_t \parallel 
                    +  \chi_{K_{\mu_{+}}} (y)  + \  \chi_{ K^{*}_{0}} (x) 
\label{bipofric}
\end{equation}
Equivalently, the relation  $ y \in \partial b_{\mu} (\cdot, y)$, 
$\mu \in \left[ \mu_{-}, \mu_{+} \right]$, with $b_{\mu}$ defined by
(\ref{Coulomb bipotential}),  
can be put in the form $y \in \partial b(\cdot,y)(x)$, with $b$ the bipotential 
defined by (\ref{bipofric}).
\label{pfric}
\end{prop}

\paragraph{Proof.} 
The collection $\left\{ b_{\mu} \mbox{ : } \mu \in
[\mu_{-}, \mu_{+}] \right\}$ of  bipotentials $b_{\mu}$ defined as in 
 (\ref{Coulomb bipotential}) is  a bipotential convex cover. 
The demonstration of the conditions of Definition \ref{defcover} is similar to 
the one of the blurred plasticity law and is not reproduced here. 

Applying Theorem \ref{thm2}, the graph (\ref{blurred Coulomb friction contact
law}) of the blurred law admits the 
bipotential:
$$b(x,y)  \ = \ \inf \left\{ \mu y_{n} \parallel x_{t} \parallel + \chi_{K_{\mu}} (y) \mbox{ : } \mu \in [\mu_{-}, \mu_{+}] \right\}\ + \  \chi_{ K^{*}_{0}} (x) \ .$$
We shall compute this bipotential. 
If $x\in K^{*}_{0}$ , three events have to be considered:
\begin{enumerate}
 \item[-] for $\parallel y_t \parallel >  \mu_{+} y_n $, $ b (x,y) \ = +\infty\ $ , 
 \item[-] for $\mu_{-} y_n \leq \parallel y_t \parallel \leq  \mu_{+} y_n $, 
                    $b(x,y)  \ = \ \inf \left\{ \mu y_n \|x_t\|  \mbox{ : } 
                                    \parallel y_t \parallel < \mu y_n  \leq  \mu_{+} y_n \right\} 
                                     \ = \ \|y_t\| \|x_t\| \ $, 
 \item[-] for $\parallel y_t \parallel \leq  \mu_{-} y_n $, 
                    $b(x,y)  \ = \ \inf \left\{ \mu y_n \|x_t\|  \mbox{ : }   
                              \mu_{-} \leq \mu \leq  \mu_{+} \right\}  
                             \ = \  \mu_{-} y_n \|x_t\|  \ $.
\end{enumerate}
If $x$ does not belongs to  $K^{*}_{0}$ , $ b (x,y) \ = +\infty\ $.  
In short, the blurred Coulomb friction contact law given by the graph
\ref{blurred Coulomb friction contact law} admits the bipotential 
(\ref{bipofric}). \quad $\square$
\section{Conclusion} 

We showed how the bipotential method is able to capture blurred 
constitutive laws.  The blurring of a constitutive law enters by   allowing  
some error in the  law with a fixed margin of tolerance. 
For example a material parameter takes indeterminate values in a 
compact set modeling some experimental incertitudes. 
The notion of bipotential convex cover leads to the construction of a 
bipotential for such a blurred constitutive law. These bipotentials 
 could be used to represent structural mechanics problems with 
 blurred constitutive laws by means of variational inequalities.

\vspace{\baselineskip}


\begin{thebibliography}{10}


\bibitem{bodo} G. Bodovill\'e: On damage and implicit standard materials,  C. R. Acad.  Sci., 
 Paris,  S\'er. II, Fasc. b, M\'ec. Phys. Astron. 
 327(8) (1999) 715-720.

\bibitem{bodo sax EJM 01} G. Bodovill\'e, G. de Saxc\'e: Plasticity with non linear kinematic 
hardening : modelling and shakedown analysis by 
the bipotential approach, Eur. J. Mech., A/Solids,  20 (2001) 99-112.


\bibitem{boush chaa IJMS 02} L. Bousshine, A. Chaaba, G. de Saxc\'e: Plastic limit load of plane 
frames with frictional contact supports,
 Int. J. Mech. Sci.  44(11) (2002)  2189-2216.


\bibitem{bipo1} M. Buliga, G. de Saxc\'e, C. Vall\'ee: Existence and construction  of 
bipotentials for graphs of  multivalued laws,   J.  Convex 
Analysis 15(1) (2008) 87-104. 

\bibitem{bipo2} M. Buliga, G. de Saxc\'e, C. Vall\'ee: Bipotentials for non 
monotone multivalued operators: fundamental results and applications, 
Acta Applicandae Mathematicae (2009), DOI 10.1007/s10440-009-9488-3. 

\bibitem{bipo3} M. Buliga, G. de Saxc\'e, C. Vall\'ee: Non maximal cyclically 
monotone graphs and construction of a bipotential for the Coulomb's 
dry friction law,   J.  Convex 
Analysis 17(1) (2010). 

\bibitem{feng hjiaj CM 06} Z.-Q. Feng, M. Hjiaj, G. de Saxc\'e, Z. Mr\'oz: Effect of frictional 
anisotropy on the quasistatic 
motion of a deformable solid sliding on a planar surface, Comput. Mech. 37 (2006) 349-361.



\bibitem{fort hjiaj CG 02} J. Fortin, M. Hjiaj, G. de Saxc\'e: An improved discrete element 
method based on a variational formulation of the 
frictional contact law, Comput. Geotech. 29(8) (2002) 609-640.

\bibitem{Halp JM 75} B. Halphen, Nguyen Quoc Son: Sur les mat\'eriaux standard g\'en\'eralis\'es, 
J. M\'ec., Paris 14 (1975) 39-63.

 

\bibitem{hjiaj bodo CRAS 00} M. Hjiaj, G. Bodovill\'e, G. de Saxc\'e: Mat\'eriaux viscoplastiques et loi 
de normalit\'e implicites, C. R. Acad.  Sci.,  Paris,  S\'er. II, Fasc. b, M\'ec. Phys. 
Astron. 328 (2000) 519-524.


\bibitem{hjiaj feng IJNME 04} M. Hjiaj, Z.-Q. Feng, G. de Saxc\'e, Z. Mr\'oz: Three dimensional finite 
element computations for frictional contact 
problems with on-associated sliding rule, Int. J. Numer. Methods  Eng. 60(12) (2004) 2045-2076.


\bibitem{laborde} P. Laborde, Y. Renard: Fixed points strategies for elastostatic frictional contact 
problems. Math. Meth. Appl. Sci. 31 (2008) 415-441.

\bibitem{saxfeng} G. de Saxc\'e, Z.Q. Feng: New inequation and functional for contact with friction: 
the implicit standard material approach, 
Mech.  Struct. and Mach. 19(3) (1991) 301-325.

\bibitem{sax CRAS 92} G. de Saxc\'e: Une g\'en\'eralisation de l'in\'egalit\'e de Fenchel et ses applications 
aux lois constitutives, 
C. R. Acad.  Sci.,  Paris,  S\'er. II 314 (1992) 125-129. 

\bibitem{sax boush KIELCE 93} G. de Saxc\'e, L. Bousshine: On the extension of limit analysis theorems to the 
non associated flow rules in 
soils and to the contact with Coulomb's friction, in: Proc. XI Polish Conference on Computer Methods in 
Mechanics  (Kielce, 1993), Vol. 2 (1993)
 815-822.

\bibitem{sax BOSTON 95} G. de Saxc\'e: The bipotential method, a new variational and numerical treatment of 
the dissipative laws of materials, in: 
Proc. 10th Int. Conf. on Mathematical and Computer Modelling and Scientific Computing, (Boston, 1995). 


\bibitem{sax boussh 2} G. de Saxc\'e, L. Bousshine: Implicit standard materials, in: Inelastic behaviour of structures under variable repeated loads, D. Weichert G. Maier (eds.), CISM  Courses and Lectures 432, Springer, Wien (2002).

\bibitem{sax feng IJMCM 98} G. de Saxc\'e, Z.-Q. Feng: The bipotential method: a constructive approach to design the complete contact law with friction and improved numerical algorithms, Math. Comput.  28(4-8) (1998) 225-245.

\bibitem{vall leri CONST 05} C. Vall\'ee, C. Lerintiu, D. Fortun\'e, M. Ban, G. de Saxc\'e: Hill's bipotential, in: New Trends in Continuum 
Mechanics, M. Mihailescu-Suliciu (ed.), Theta Series in Advanced Mathematics, Theta Foundation, Bucarest (2005) 339-351.

\bibitem{ER1} I. Babuska, W.C. Rheinboldt: A posteriori estimates for the finite element method, Int. J. Numer. Methods Engrg. 12 (1978) 1597-1615.

\bibitem{ER2} I. Babuska, W.C. Rheinboldt: Adaptive approaches and reliability estimators in finite element analysis, Comput. Methods Applied Mech. Engrg. 17/18 (1979) 519-540.

\bibitem{ER3} D.W. Kelly, J. Gago, O.C. Zienkiewicz, I. Babuska, A posteriori error analysis and adaptative processes in finite element method. Part 1: Error analysis, Int. J. Numer. Methods Engrg. 19 (1983) 1593-1619.

\bibitem{ER4} J. Gago, D.W. Kelly, O.C. Zienkiewicz, I. Babuska, A posteriori error analysis and adaptative processes in finite element method. Part 2: Adaptative mesh refinement, Int. J. Numer. Methods Engrg. 19 (1983) 1921-1656.

\bibitem{ER5} J.T. Oden, L. Demkowicz, W. Rachowicz, T.A. Westermann, Toward a universal h-p adaptative finite element strategy, Part 2: A posteriori error estimation, Comput. Methods Appl. Mech. Engrg. 77 (1989) 113-180.

\bibitem{ST1} O.C. Zienkiewicz, J.Z. Zhu, A simple error estimator and adaptative procedure for practical engineering analysis, Int. J. Numer. Methods Engrg. 24 (1987) 337-357.

\bibitem{ST2} O.C. Zienkiewicz, J.Z. Zhu, The superconvergent patch recovery and a posteriori error estimates, Part 1: The recovery technique, Int. J. Numer. Methods Engrg. 33 (1992) 1331-1364.

\bibitem{ST3} O.C. Zienkiewicz, J.Z. Zhu, The superconvergent patch recovery and a posteriori error estimates, Part 2: Error estimates and adaptativity, Int. J. Numer. Methods Engrg. 33 (1992) 1365-1382.

\bibitem{DA1} B. Fraeijs de Veubeke, Displacement and equilibrium models in the finite element method, in: O.C. Zienkiewicz, ed., Stress Analysis Chap. 9 (John Wiley, 1965).

\bibitem{DA2} J.F. Debongnie, H.G. Zhong, P. Beckers, Dual analysis with general boundary conditions, Comput. Methods. Appl. Mech. Engrg. 122 (1995) 183-192.

\bibitem{LAD1} P. Ladev\`eze, Comparaisons de mod\`eles de milieux continus, Th\`ese d'Etat, Universit\'e Pierre et Marie Curie, Paris (1975).

\bibitem{LAD2} P. Ladev\`eze, G. Coffignal, J.P. Pelle, Accuracy of elastoplastic and dynamic analysis, in: I. Babuska, J. Gago, E. Oliveira and O.C. Zienkiewicz, eds., Accuracy Estimates and Adaptative Refinements in Finite Element Computations (John Wiley, 1986) 181-203.

\bibitem{LAD3} P. Ladev\`eze, J.P. Pelle, P. Rougeot, Error estimation and mesh optimization for classical finite element, Engrg. Comput. 8 (1991) 69-80.

\bibitem{LAD4} P. Ladev\`eze, La ma\^itrise des mod\`eles en m\'ecanique des structures, Revue europ\'eenne des \'el\'ements finis 1 (1992) 9-30.

\bibitem{LAD5} P. Ladev\`eze, N. Moes, A new a posteriori error estimation for nonlinear time-dependent finite element analysis, Comput. Methods. Appl. Mech. Engrg.  157 (1997) 45-68.

\bibitem{LAD6} P. Ladev\`eze, J.P. Pelle, La ma\^itrise du calcul en m\'ecanique lin\'eaire et non lin\'eaire, (Hermes Science, 2001).

\bibitem{LAD7} Pierre Ladev\`eze, E. Florentin, Verification of stochastic models in uncertain environments
using the constitutive relation error method, Comput. Methods. Appl. Mech. Engrg.  196 (2006) 225-224.

\bibitem{LAD8} P. Ladev\`eze, G. Puel, A. Deraemaeker, T. Romeuf, Validation of structural dynamics models containing uncertainties, Comput. Methods. Appl. Mech. Engrg.  195 (2006) 373-393.

\end{thebibliography}
\end{document}